\documentclass[a4paper,12pt]{amsart}

\usepackage{amssymb}

\usepackage{amsfonts}

\def\C{\mathbb{C}}
\def\Z{\mathbb{Z}}
\def\M{\mathcal{M}}
\def\V{\mathcal{V}}
\def\d{\partial}
\def\oM{\overline{\mathcal{M}}}

\usepackage{epsf}
\newcommand{\inspic}[1]{\begin{tabular}{c}\epsfbox{#1}\end{tabular}}

\newtheorem{proposition}{Proposition}

\title[BCOV theory via Givental group action on CohFT]{BCOV theory via Givental group action on cohomological fields theories}

\author{Sergey Shadrin}

\address{Korteweg-de~Vries Institute for Mathematics, University of Amsterdam,
Plantage Muidergracht 24, 1018 TV Amsterdam, The Netherlands}

\email{s.shadrin@uva.nl}

\address{Department of Mathematics, Institute of System Research, Nakhimovsky prospekt 36-1, Moscow 117218, Russia}

\email{shadrin@mccme.ru}

\begin{document}

\begin{abstract}
In a previous paper, Losev, me, and Shneiberg constructed a full descendant
potential associated to an arbitrary cyclic Hodge dGBV algebra. This contruction
extended the construction of Barannikov and Kontsevich of solution of the WDVV
equation, based on the earlier paper of Ber\-shad\-sky, Ce\-cot\-ti, Oo\-gu\-ri, and Va\-fa.

In the present paper, we give an interpretation of this full descendant potential 
in terms of Givental group action on cohomological field theories. In particular, the
fact that it satisfies all tautological equations becomes a trivial observation.
\end{abstract}

\maketitle

\tableofcontents

\section{Introduction}

In this paper, which is partly a survey and introductory paper, we revise our
previous work~\cite{LosShaShn07} on Hodge field theory in terms of the Givental
group action on
cohomological fields theories. Let us first remind the history and motivation of
the whole subject.

\subsection{BCOV theory}

In~\cite{BerCecOogVaf94}, Bershadsky, Cecotti, Ooguri, and Vafa introduced an
approach to construct a mirror partner for the Gromov-Witten potential of a
mirror dual Calabi-Yau manifold. They inroduced what is called now BCOV action
on the space of homorphic polyvector fields with coefficients in antiholomorphic
forms on a Calabi-Yau manifold.

In~\cite{BarKon98}, Barannikov and Kontsevich proved that indeed the critical
value of BCOV action is a solution of the WDVV equation. In fact, their argument
works in the settings of an abstract Hodge dGBV algebra (the algebraic
structure that formalizes the properties of the space of polyvector fields on
Calabi-Yau manifolds), see~\cite{Man99a, Man99b}.

In~\cite{LosSha07}, Losev and me represented the Barannikov-Kontsevich
solution of the WDVV equation as a sum over trivalent trees and gave a new proof
that it satisfies the WDVV equation. Also we considered its genus expansion 
(the perturbative expansion of the BCOV action). It has
appeared that the additional $1/12$-axiom allows to prove some other PDEs coming from the geometry of the moduli spaces
of curves (in addition to
WDVV) for the higher genus amplitudes. 

In order to explain this additional axiom we introduced~\cite{LosSha07} a special TCFT on the
real blow-up of the moduli space of curves (Kimura-Stasheff-Voronov
space~\cite{KimStaVor95}). Following ideas of Getlzer, we showed~\cite{LosSha07}
that indeed, the relations between the Dehn twists in the mapping class group
imply the properties of the BV operator in dGBV algebra and $1/12$ axiom in the
same fashion~\cite{Get94}; and using the induction procedure closed to Zwiebach's idea of
splitting the moduli space into the interior part and the boundary~\cite{Zwi93}, we got a
geometric construction that could degenerate to the expansion of BCOV action in
the simplest possible case.

This gave a hint of how one can introduce descendants in order to extend BCOV
theory to the full descendant potential. First, one could do this by hand in
genus $0$ using the topological recursion relation, second one could introduce descendants in this special TCFT (called `Zwiebach
invariants' in~\cite{LosSha07, LosShaShn07}) and look what happens at
degeneration. Both approaches give the same answer (in terms of graphs), 
and after a sequence of extrimely complicated calculations with
graphs~\cite{Sha07,ShaShn07,Shn08} Losev, me, and Shneiberg proved
in~\cite{LosShaShn07} that the full descendant potential that we defined indeed
satisfies all expected relations coming from geometry of the moduli space of
curves.

The complete construction is called in~\cite{LosShaShn07} Hodge field theory.

\subsection{Givental group action}

Givental group action on cohomological field theories has first appeared in the papers
of Givental~\cite{Giv01a,Giv01,Giv04}. He observed that the localization formulas for
the Gromov-Witten potentials of projective spaces and some Fano complete
intersections have the structure that can be formulated in terms of the
quantization of some group action on genus $0$ parts of potentials. 

He explained the group action in genus $0$ in terms of symplectic geometry. In
particular, all semi-simple theories form one orbit of this action (or a small
family of orbits~\cite{Kaz07}). Then, using the quantization procedure for the
group action, Givental suggested a universal formula for the genus expansion of
an arbitrary Frobenius manifold.

In particular, Givental conjectured that his universal formula, when applied to
the genus $0$ Gromov-Witten potential of a target manifold with the semi-simple
quantum cohomology ring, reconstructs the full Gromov-Witten potential. This
conjecture was proved in~\cite{Tel07} via a complete classification of
semi-simple cohomological field theories.

Meanwhile, Faber, me, and Zvonkine proved in~\cite{FabShaZvo06} that after
quantization Givental's group action preserves all universal constrains coming
from geometry of the moduli space of curves (tautological relations). In
particular, if one applies a (quantized) element of the Givental group to 
the potential of a cohomological field theory, one gets the potential of another
cohomological field theory.

Our argument was generalized by Kazarian~\cite{Kaz07}, who explained a part of
the Givental group action as an action on the cohomology classes in the moduli
space of curves.

\subsection{Plan of the paper}

In Section~\ref{sec:notions}, we briefly introduce all necessary basic notions such as
moduli space of curves, cohomological field theory, tautological relations and
so on. In Section~\ref{sec:givental}, we discuss different ways to introduce a
part of Givental's action (this part is called `upper triangular group' in the
literature). In Section~\ref{sec:oldhodge}, we remind the construction of the
full descendant potential of Hodge field theory. Finally, in
Section~\ref{sec:newhodge}, we give an alternative construction of this
potential using the Givental group action.

\subsection{Remarks}

There are several remarks. First, there is a particular mathematical trouble --
usually, the cyclic Hodge dGBV algebras are infinite-dimentional, and there is a
problem of convegence of all the tensor expressions used in
Sections~\ref{sec:oldhodge} and~\ref{sec:newhodge}. We don't worry about it, since
there exist different approaches to renormalization of BCOV
theory~\cite{Cos07a,Dij95}, and they
seem to be very well compatible with our technique.

Second, we don't know, whether $1/12$-axiom is valid in known examples of Hodge
dGBV algebras. It has never appeared before since people used to study only
genus $0$. In the only example~\cite{Dij95}, where all computations can be done explicitely
(elliptic curve), it is valid, but in appears to be a statement of the type
$0=0$, so the coefficients are not essential. 

Third, there is another family of examples of Hodge dGBV algebras. It appears to
be the structure on the space of differential forms on symplectic manifolds
satisfying the hard Lefschetz condition~\cite{Mer98}. However, these examples
seems not be studied yet elsewhere.

Forth, the axioms of Hodge dGBV algebra also have another motivation. Losev's
description~\cite{Los98} of Frobenius structures on the base space of the universal 
unfolding of a simple singularity leads to some algebraic formalism, whose
immediate non-linear generalization appears to be exactly the axioms of Hodge dGBV
algebra.

Fifth, the full descendant potential that we study here (a genus expansion of
Barannikov-Kontsevich construction) is not exactly the mirror partner of the
Gromov-Witten potential of a mirror dual Calabi-Yau manifold. Barannikov
observed~\cite{Bar01,Bar02} that there is a family of possible changes of variables, 
parametrized by semi-infinite variations of Hodge structure, and one of them is required to
complete the mirror construction. In fact, Barannikov construction seems to be
nothing but again the upper-triangular part of the Givental group action with
the additional restriction that is preserves the weights of homogeneity of the
formal variables in the potential. We are going to discuss it elsewhere.

\subsection{Acknowledgements}
I am very grateful to E.~Feigin, E.~Getzler, M.~Kazarian, A.~Losev, I.~Shneiberg, and
D.~Zvonkine for the plenty of fruitful discussions. The remarks of the anonymous referee were very helpful.

\section{Basic notions}\label{sec:notions}

\subsection{The moduli space of curves}

The moduli space of curves $\M_{g,n}$, $g\geq 0$, $n\geq 0$, $2g-2+n>0$, parametrizes smooth 
complex curves of genus $g$ with $n$ ordered marked points. It is a smooth complex orbifold of dimension $3g-3+n$. 

The space $\oM_{g,n}$ is a compactification of $\oM_{g,n}$. It parametrizes stable curves of 
genus $g$ with $n$ ordered marked points. A stable curve is a possibly reducible curve with possible 
nodes, such that the order of its automorphism group is finite. Genus of a stable curve is the arithmetic 
genus, namely, the genus of the smooth curve that we get if the replace each node (given locally by the 
equation $xy=0$) with a cylinder (given locally by the equation $xy=\epsilon$). The space $\oM_{g,n}$ 
is a smooth compact complex orbifold.

The space $\oM_{g,n}$ has a natural stratification by the topological type of stable curves. 

\subsection{Natural mappings}

There is a number of natural mappings between the moduli spaces of curves.

First, there are projections $\pi\colon \oM_{g,n+1}\to \oM_{g,n}$ that forget the last marked point. 
Note that there is a subtlety related to the fact that when we forget a marked point a stable curve 
can become unstable.

Second,
there is a 2-to-1 mapping $\sigma\colon\oM_{g-1,n+2}\to\oM_{g,n}$ whose image it the boundary
divisor of irreducible curves with one node. Third, there are mappings
$\rho\colon \oM_{g_1,n_1+1}\times \oM_{g_2,n_2+1}\to \oM_{g,n}$, $g_1+g_2=g$,
$n_1+n_2=n$, whose images are the other irreducible boundary divisors of the
compactification of $\oM_{g,n}$.

That gives a complete description of strata in codimension 1. Also we see that
the maps $\rho_*$ and $\sigma_*$ induce the structure of modular
operad~\cite{GetKap98} on the (co)homologies of the moduli spaces of curves.

\subsection{Strata in $\oM_{g,n}$}\label{sec:strata}

Any irreducible boundary stratum of codimension $k$ in $\oM_{g,n}$ is represented 
as an image $p(S)$ of the product $S=\oM_{g_1,n_1}\times \cdots \times \oM_{g_a, n_a}$. 
Here $p$ is a composition of the $k$ mappings $\sigma$ and/or $\rho$ described above. 
We can associate a dual graph $G$ to the mapping $p$. The vertices $v_i$, $i=1,\dots,a$ 
correspond to the spaces $\oM_{g_i,n_i}$; $v_i$ is marked by $g_i$ (a non-negative integer 
that keeps track of genus) and its index is $n_i$. The marked points on curves in $\oM_{g_i,n_i}$ 
that remain marked point in the image of $p$ correspond to the leaves (free half-edges) of $G$. 
The pairs of marked points that turn into nodes in the image of $p$ form edges of the graph $G$.

We see that there are exactly $k$ edges in $G$ and $G$ is connected. Moreover, the 
requirements on the arithmetic genus and the number of marked points of curves in 
the image of $p$ mean that $\sum_{i=1}^a g_i + b_1(G)=g$ and $\sum_{i=1}^a n_i =n+2k$.

\subsection{Tautological classes}

The cohomology of the moduli space of curves is
a complicated object that is not yet studied well enough. Instead of it, people
usually consider a special system of subalgebras of the cohomology algebras
of the moduli spaces called tautological rings.

The system of tautological rings 
$RH^*(\oM_{g,n})\subset H^*(\oM_{g,n})$
is defined as a minimal system of subalgebras of the algebras of cohomology that
is closed under the
push-forwards and pull-backs via the natural mappings between the moduli spaces.

The cohomology classes in $RH^*(\oM_{g,n})$ are called tautological classes. It
follows immediately from the definition that $RH^*(\oM_{g,n})$ also form a
modular operad.

\subsection{Additive generators of the tautological ring}

We describe a system of additive generators of the tautological ring of the space $\oM_{g,n}$. 

Let $L_i$ the line bundle over $\oM_{g,n}$, whose fiber over a point $x\in
\oM_{g,n}$ represented by a curve $C_g$ with
marked points $x_1,\dots,x_n$ is equal to $T^*_{x_i}C_g$.
Denote by $\psi_i\in H^2(\oM_{g,n})$ the first Chern class of $L_i$. It is easy
to show that $\psi_i$ is a tautological class.
Denote by $\kappa_j$, $j=1,2,\dots,$ the class $\pi_*(\psi_{n+1}^{j+1})$,
$\pi\colon \oM_{g,n+1}\to \oM_{g,n}$. Since $\psi$-classes are tautological, 
$\kappa_j\in RH^{2j}(\oM_{g,n})$.

Now the system of additive generators of $RH^*(\oM_{g,n})$ can be described in
terms of $\psi$-$\kappa$-strata. Let $p(S)$ be an irreducible stratum of
codimension $k$ in $\oM_{g,n}$, $S=\oM_{g_1,n_1}\times \cdots \times \oM_{g_a,
n_a}$, as in Section~\ref{sec:strata}. Equip each $\oM_{g_i,
n_i}$, $i=1,\dots,a$, with a monomial of $\psi$- and $\kappa$-classes. Consider
push-forward $p_*$ of the class that we get in the cohomology of $S$. 
What we
obtain is called $\psi$-$\kappa$-stratum, and all
$\psi$-$\kappa$-strata form a system of additive generators of
$RH^*(\oM_{g,n})$.

\subsection{Tautological relations}\label{sec:tau-rel}

Of course, $\psi$-$\kappa$-strata are not free additive generators. There are
plenty of relations that are called tautological relations. Let us give a few
examples. In examples, we use dual graphs defined in Section~\ref{sec:strata},
decorated by $\psi$- and $\kappa$-classes.

In $\oM_{0,4}$ we have: $\psi_1=\psi_2=\psi_3=\psi_4=\kappa_1$. Moreover all
these classes are also equal to
\begin{equation}
\inspic{pictures.1}=\inspic{pictures.2}=\inspic{pictures.3}.
\end{equation}
Leaves of graphs are labeled by the numbers of the corresponding marked points.

In $\oM_{1,1}$ we have:
\begin{equation}
\inspic{pictures.4}=\inspic{pictures.5}=\frac{1}{12}\inspic{pictures.6}.
\end{equation}

There is a very small number of basic tautological relations known
explicitely~\cite{Wit01,Kee92,Mum83,Get97,Get98,BelPan00,KimLiu06}. 
Also, there are general theorems proving the existence of some special families of
tautological relations~\cite{Get98,Ion02,FabPan05}.

\subsection{Gromov-Witten theory}

Gromov-Witten theory associates to a compact K\"ahler manifold $X$ some algebraic
structure on its cohomology $H^*(X)$ that generalizes the standard structure
of a graded Frobenius algebra (that is, a structure of graded commutative
associative algebra with a unit compatible with a non-degenerate even scalar
product). From the algebraic point of view, this structure on a vector space $V$
($=H^*(X)$)
equipped with a scalar product is given by a unital representation of a modular
operad of the
cohomologies $H^*(\oM_{g,n})$ of the moduli spaces of curves $\oM_{g,n}$. 

Let us describe the operations on $H^*(X)$. We consider
\begin{equation}
\alpha_{g,n}^{u,\beta}(v_1,\dots,v_n):=\int_{[\oM_{g,n}(X,\beta)]^{vir}}f^*(u)ev_1^*(v_1)\cdots
ev_1^*(v_1).
\end{equation}
as a linear map $H^*(X)^{\otimes n}\to \C$; it depends on $g\geq 0$, $\beta\in
H_2(X)$, and $u\in H^*(\oM_{g,n})$.

Here $\oM_{g,n}(X,\beta)$ is the moduli space of stable maps of curves of genus
$g$ with $n$ marked points to $X$, such that the image of the fundamental class
is equal to $\beta\in H_2(X)$. Such spaces consist of several irreducible
components of different dimension, so the integral is taken over a particular
homology class $\oM_{g,n}(X,\beta)$ that is calles the virtual
fundamental class and replace the usual fundamental class.

There is a forgetful map $f\colon \oM_{g,n}(X,\beta)\to\oM_{g,n}$ that associate
to a stable map $(C_g,x_1,\dots,x_n,\phi\colon C_g\to X)$ the stabilization of
the source curve with marked points. The map $ev_i\colon \oM_{g,n}(X,\beta) \to
X$, $i=1,\dots, n$ is defined as $\phi(x_i)$.

\subsection{CohFT}

The main properties of Gromov-Witten theory are captured by the notion of
cohomological field theory (CohFT). Roughly speaking, a CohFT is a system of
cohomology classes on the moduli spaces of curves with the values in the tensor 
powers of $V$, compatible with all natural mappings between the moduli spaces. 

The formal definition is the following. We fix a vector space $V=\langle
e_1,\dots, e_s\rangle$ ($e_1$ will play a special role) with a
non-degenerate scalar product $\eta$. A cohomological field theory is a system of
cohomology classes $\alpha_{g,n}\in H^*(\oM_{g,n},V^{\otimes n})$ satisfying the
properties:
\begin{enumerate}
\item $\alpha_{g,n}$ is equivariant with respect to the action of $S_n$
on the labels of marked point and component of $V^{\otimes n}$.
\item $\sigma^*\alpha_{g,n}=\left(\alpha_{g-1,n+2},\eta^{-1}\right)$;
$\rho^*\alpha_{g,n}=\left(\alpha_{g_1,n_1+1}\cdot\alpha_{g_2,n_2+1},\eta^{-1}\right)$
 (in both cases we contract with the scalar product the two components of $V$ 
 corresponding to the two points in the preimage of the node under
 normalization).
\item 
$\left(\alpha_{0,3},e_1\otimes e_i\otimes
e_j\right)=\eta_{ij}$, 
$\pi^*\alpha_{g,n}=\left(\alpha_{g,n+1}, e_1\right)$
(again, we contract the component of $V$ corresponding to the last marked
point with $e_1$).
\end{enumerate}

\subsubsection{TFT}\label{sec:TFT}

Topological field theory is a special case of CohFT, when all cohomology classes
$\alpha_{g,n}$ are in degree $0$. One can easily show that in this case the whole
system $\{\alpha_{g,n}\}$ is determined by $\alpha_{0,3}\in V^{\otimes 3}$. Of
course, there are some restrictions for the choice of the $3$-tensor $\alpha_{0,3}$. One
can show that a particular $3$-tensor can be exteded to a TFT if an only if this
$3$-tensor determines a structure of commutative Frobenius algebra with the unit
on the vector space $V$ with the scalar product $\eta$~\cite{Dub96}.

\subsection{Correlators}\label{sec:correlators}

We want to assign to a CohFT a numerical information that would allow to reconstruct 
completely its restriction to the tautological ring. Actually, it appears that
it is enough to know all the integrals of the type
\begin{equation}
\langle \tau_{d_1}(e_{i_1})\cdots\tau_{d_n}(e_{i_n}) \rangle_g 
:=\int_{\oM_{g,n}}\prod_{j=1}^n\psi_j^{d_j}\cdot \left(\alpha_{g,n},\otimes_{j=1}^ne_{i_j}\right)
\end{equation}
(correlators) in order to be able to reconstruct
$\int_{\oM_{g,n}}\beta \cdot \left(\alpha_{g,n},\otimes_{j=1}^ne_{i_j}\right)$ for an
arbitrary tautological class $\beta$.

Indeed, using the property (3) of CohFT, we can express the integrals with the
monomials of $\psi$- and $\kappa$-classes via the integrals with the monomials
of $\psi$-classes only. Then, using the property (2), we can combine the
integrals with monomials of $\psi$- and $\kappa$-classes in order to get the
integrals with the class $\beta$ represented with an arbitrary dual graph
decorated by $\psi$- and $\kappa$-classes.

Usually, it is convenient to gather the correlators into a generating series
(potential of CohFT):
\begin{align}
Z & =\exp\left(\sum_{g\geq 0} \hbar^{g-1}F_g\right) \\
 & =\exp\left(\sum_{g,n} \hbar^{g-1} \sum \frac{1}{n!}
\sum_{\mathbf{d},\mathbf{i}}\langle \prod_{j=1}^n \tau_{d_j}(e_{i_j})
\rangle_g \prod_{j=1}^n t_{d_j,i_j} \right). \notag
\end{align}
Here $t_{d,i}$ are some formal variables. In the second line, the first sum is taken over $g,n\geq 0$, $2g-2+n>0$,
and the second sum is taked over all possible tuples of indices
$\mathbf{d}=(d_1,\dots,d_n)$, $\mathbf{i}=(i_1,\dots,i_n)$.

\subsection{Tautological equations}\label{sec:tauteq}

First, there is a freedom in the choice of expressions of the intergrals with an arbitrary 
tautological class in terms of correlators. Second, $\psi$-$\kappa$-strata are not free 
additive generators of the tautological ring, as discussed in Section~\ref{sec:tau-rel}. This implies 
that correlators should satisfy a huge number of universal relations. `Universal' means that 
the realtions depend only on the choice of the vector space with the scalar product and the unit.

These universal relations can be gathered in an infinite number of very complicated PDEs on
 the series $F_g$, $g\geq 0$. These PDEs are the following:
\begin{align}
& \exp(-\sum t_{0,i}\eta^{ij} t_{0,j})\frac{\d Z}{\d t_{0,1}} = \sum t_{d+1,i}\frac{\d Z}{\d t_{d,i}} & (\mathrm{string}); \\
& \frac{\d F_g}{\d t_{1,1}} - (2g-2)F_g = \sum t_{d,i}\frac{\d F_g}{\d t_{d,i}}, \quad g\geq 0 & (\mathrm{dilaton});
\end{align}
and a system of PDEs associated to each tautological relation in the sense of Section~\ref{sec:tau-rel}. 
Each system of PDEs collects the relations for correlators that are implied by all pull-backs of a given 
relation multiplied by an arbitrary monomial of $\psi$-classes. For example, the relation
\begin{equation}
\inspic{pictures.7}=\inspic{pictures.1}
\end{equation}
implies the following system of PDEs for $F_0$:
\begin{equation}
\frac{\d^3 F_0}{\d t_{a+1,i}\d t_{b,j}\d t_{c,k}}=\sum \frac{\d^2 F_0}{\d t_{a,i}\d t_{0,\alpha}}\eta_{\alpha,\beta}
\frac{\d^3 F_0}{\d t_{0,\beta}\d t_{b,j}\d t_{c,k}}
\end{equation}
(topological recursion relation in genus 0, or just TRR-0 for short).


\section{Givental's group action on CohFT}\label{sec:givental}

In this section we describe a part of the Givental theory of the group action on
CohFT~\cite{Giv01,Giv04}. 

\subsection{Action on cohomology classes}

Let $R(z)=Id+zR_1+z^2R_2+
\dots \in Hom(V,V)\otimes\C[[z]]$ be a series of endomorphisms of the space $V$ 
satisfying $R^*(-z)R(z)=Id$. In other words, $R(z)=\exp(r_1z+r_2z^2+\dots)$, 
where $r_l\in Hom(V,V)$ is (graded) symmetric matrix for odd $l$ and skewsymmetric 
for even $l$, with respect to the scalar product $\eta$.

Following~\cite{Kaz07}, we associate to $\sum_{l=1}^{\infty}r_lz^l$ an infinitesimal deformation of CohFT. 
Denote by $(r_lz^l)\hat{\ } \alpha_{g,n}$ the following class on $\oM_{g,n}$:
\begin{align}
& -\pi_*\left(\alpha_{g,n+1}\psi_{n+1}^{l+1},r_l(e_1)\right)
+\sum_{k=1}^n\psi_k^l
r_l^{(k)}\alpha_{g,n} \label{eq:kazarian}  \\
&
+\frac{1}{2}\sum_{i=0}^{l-1}(-1)^{i+1}\sigma_*\left(\alpha_{g-1,n+2}\psi_{n+1}^i\psi_{n+2}^{l-1-i},
\eta^{-1}r_l\right)
\notag \\
& +\frac{1}{2}\sum_{div}\sum_{i=0}^{l-1}(-1)^{i+1}
\rho_*\left(\alpha_{g_1,n_1+1}\psi_{n_1+1}^i\cdot
\alpha_{g_2,n_2+1}\psi_{n_2+1}^{l-1-i},\eta^{-1}r_l\right).\notag
\end{align}
The last sum here ($\sum_{div}$) is taken over all irreducible boundary divisors, whose generic
points are represented by two-component curves.

We define 
$\alpha'_{g,n}:=\exp\left(\sum_{l=1}^\infty(r_lz^l)\hat{\ }\right)\alpha_{g,n}$
applying successively formula~\eqref{eq:kazarian} simultaneously to the whole system of
classes, for all $g$ and $n$.
Kazarian~\cite{Kaz07} proved that (1) $\alpha'_{g,n}$ are well-defined
cohomology classes with the values in the tensor powers of $V$ and (2)
$\alpha'_{g,n}$ form a CohFT.

\subsection{Action on correlators}\label{sec:action-corr}

We start with the same operator $R(z)=\exp(r_1z+r_2z^2+\dots)$ as before,
but now we denote by $(r_lz^l)\hat{\ }$ the following differential operator:
\begin{align}
& -(r_l)^1_\mu\frac{\d}{\d t_{l+1,\mu}}
+ \sum_{d=0}^\infty t_{d,\nu} (r_l)^\nu_\mu \frac{\d}{\d t_{d+l,\mu}}
\label{eq:derivative} \\
& +\frac{\hbar}{2} \sum_{i=0}^{l-1}(-1)^{i+1} (r_l)_{\mu,\nu}\frac{\d^2}{\d t_{i,\mu}\d
t_{l-1-i,\nu}}. \notag
\end{align}

Let $Z=\exp\left(\sum_{g=0}^\infty \hbar^{g-1}F_g\right)$ be the potential associated to a CohFT (see
Section~\ref{sec:correlators}). Givental proved~\cite{Giv01} that 
$Z':=\exp\left(\sum_{l=1}^\infty(r_lz^l)\hat{\ }\right)Z$ is also a well-defined
formal power series of the same type as $Z$, namely $Z'$ is represented as
$\exp\left(\sum_{g=0}^\infty \hbar^{g-1}F'_g\right)$. Moreover, Faber, me, and Zvonkine
proved in~\cite{FabShaZvo06} that $Z'$ is also the potential of a CohFT.

In fact, the theorem of Kazarian mentioned above is dual to our theorem, and 
it is even more general, since it doesn't require to restrict CohFT to the
tautological ring.

\subsection{Weyl quantization}

We want to explain, following Givental, why the operator given by
Formula~\eqref{eq:derivative} (or, the dual operator on cohomology classes given
by Formula~\eqref{eq:kazarian}) is denoted by $(r_lz^l)\hat{\ }$. 

\subsubsection{Action in genus 0}

Consider the genus $0$ part $F_0$ of the potential of a CohFT. It should satisfy
only string, dilaton, and TRR-0 (see Section~\ref{sec:tauteq}).
It is possible to encode these properties in geometric terms. Let $\V$ be the 
space $V\otimes \C((z^{-1}))$ with the natural symplectic structure 
\begin{equation}
\omega(f,g)=\frac{1}{2\pi\sqrt{-1}}\oint \left(f(-z),g(z)\right)dz, \quad f(z),g(z)\in \V.
\end{equation}
This allows to identify $\V$ with the space $T^*\V_+$, $\V_+=V\otimes \C[z]$. 
The natural coordinates $q_{n,i}$ in the space $z^n V$ (dual to a 
chosen basis $\langle e_1,\dots,e_s\rangle$ in $V$) are identified with $t_{n,i}$; 
the only exception is $q_{1,1}=t_{1,1}-1$.

Givental proved~\cite{Giv01, Giv04} that the graph $\mathcal{L}$ of $dF_0$ in $\V$ is the germ of a 
Lagrangian cone with the vertex at $0$ such that its tangent spaces $L$ are 
tangent to $\mathcal{L}$ exactly along $zL$.
This immediately implies that the group of linear symplectic operators
compatible with the multiplication by $z$, that is, of the 
type $M(z)\in End(V)\otimes \C((z^{-1}))$, $M^*(-z)M(z)=Id$, 
acts on on the formal power series 
satisfying string, dilaton, and TRR-0. The group of such operators is 
called the twisted loop group. 

\subsubsection{Quantization}

Consider an element $m$ of the twisted loop group. It turns the graph of $dF_{0}$ into 
the graph of $dF'_0$, where $F'_0$ is a formal power series that satisfies
string, dilaton, and TRR-0. The logarithm $\log m$ of this element is a linear vector 
field with quadratic Hamiltonian $h_m$. We can quantize any quadratic Hamiltonian in the 
following standard way:
\begin{align}
(p_{n_1,i_1}p_{n_2,i_2})\hat{\ }&=\hbar \frac{\d^2}{\d q_{n_1,i_1}\d
q_{n_2,i_2}},  \\
(p_{n_1,i_1}q_{n_2,i_2})\hat{\ }&=q_{n_2,i_2}\frac{\d}{\d q_{n_1,i_1}}, \notag \\
(q_{n_1,i_1}q_{n_2,i_2})\hat{\ }&=\frac{1}{\hbar}q_{n_1,i_1}q_{n_2,i_2} \notag
\end{align}
(the coordinates $p_{n,i}$ are the Darboux pairs of $q_{n,i}$ with respect to the symplectic structure on $\V$).
Then one can consider the operator $\exp(\hat{h_m})$.

Applying this quantization to the operator $R(z)=Id+zR_1+z^2R_2+
\dots \in Hom(V,V)\otimes\C[[z]]$, we obtain the action on full potentials of
CohFT. Strictly speaking, it is not obvious that this action is well-defined for
the infinite series in $z$; but actually it is. The key additional property of
the full potentials of CohFT is the following: $\langle
\tau_{d_1}(e_{i_1})\cdots\tau_{d_n}(e_{i_n}) \rangle_g=0$ if $\sum_{j=1}^n
d_j>3g-3+n$. It is called $3g-2$ property~\cite{EguXio98, Get04, Giv01}. It follows from the tautological
relations $\psi_1^{d_1} \cdots \psi_n^{d_n} =0$ in $\oM_{g,n}$ if $\sum_{j=1}^n
d_j>3g-3+n$, and it can be written down as a system of PDEs. This property
guarantees that the action of $\exp(\hat{h_R})$ on the full potentials of CohFT
is well-defined.

Explicit formulas for $(r_lz^l)\hat{\ }$ were first computed and studied
in~\cite{Lee03}, see also~\cite{Lee06a,Lee06b,Lee07}.

\subsection{Remark on gradings}

Actually, below we are applying not exactly the Givental group action as it is
presented here, but a ($\Z_2$-) graded version of it. Namely, the target space
$V$ is a graded vector space. This also requires a graded version of CohFT. 
In our application these grading will affect some signes in the formulas. We
hope that the reader will be able to reconstruct the precise definitions of the
graded versions of CohFT and the Givental group action himself, while we are going to ignore
this issue.


\section{Hodge field theory}\label{sec:oldhodge}

In this section we present the algebraic part of the construction of the
potential associated to cyclic Hodge dGBV-algebra developed by Losev, me, and
Shneiberg.

\subsection{Cyclic Hodge algebras}\label{sec:cH-algebras}
In this section, we give the definition of cyclic Hodge
dGBV-alge\-bras, see~\cite{LosSha07,LosShaShn07}.
(cyclic Hodge algebras, for short). A supercommutative associative $\C$-algebra $H$ with
the unit is called cyclic Hodge algebra, 
if there are two odd linear operators $Q,G_-\colon H\to H$ and an even linear function 
$\int\colon H\to\C$ called integral. They must satisfy the following seven axioms A1-A7:

\begin{description}

\item[A1] $(H,Q,G_-)$ is a bicomplex: 
\begin{equation}
Q^2=G_-^2=QG_-+G_-Q=0;
\end{equation}

\item[A2] $H=H_0\oplus H_4$, where $QH_0=G_-H_0=0$ and $H_4$ is represented as 
a direct sum of subspaces of dimension $4$ generated by 
$e_\alpha, Qe_\alpha, G_-e_\alpha, QG_-e_\alpha$ for some vectors 
$e\in H_4$, i.~e. 
\begin{equation}
H=H_0\oplus
\bigoplus_{\alpha} 
\langle 
e_\alpha, Qe_\alpha, G_-e_\alpha, QG_-e_\alpha
\rangle
\end{equation}
(Hodge decomposition);

\item[A3] $Q$ is an operator of the first order, it satisfies the Leibniz rule: 
\begin{equation}\label{eq:leibniz}
Q(ab)=Q(a)b+(-1)^{\tilde a}aQ(b)
\end{equation}
(here and below we denote by $\tilde a$  the parity of $a\in H$);

\item[A4] $G_-$ is an operator of the second order, 
it satisfies the $7$-term relation:
\begin{align}\label{eq:7-term}
G_-(abc)& = G_-(ab)c+(-1)^{\tilde b(\tilde a+1)}bG_-(ac)
+(-1)^{\tilde a}aG_-(bc)\\
& -G_-(a)bc-(-1)^{\tilde a}aG_-(b)c
-(-1)^{\tilde a+\tilde b}abG_-(c). \notag
\end{align}

\item[A5] $G_-$ satisfies the property called $1/12$-axiom: 
\begin{equation}\label{eq:1/12-axiom}
str(G_-\circ a\cdot)=(1/12)str(G_-(a)\cdot)
\end{equation}
(here $a\cdot$ and $G_-(a)\cdot$ are the 
operators of multiplication by $a$ and $G_-(a)$ respectively, $str$ means
supertrace).

\end{description}

Define an operator $G_+\colon H\to H$. We put $G_+H_0=0$, and on each subspace 
$\langle e_\alpha, Qe_\alpha, G_-e_\alpha, QG_-e_\alpha \rangle$ 
we define $G_+$ as 
\begin{align}
& G_+e_\alpha=G_+G_-e_\alpha=0, \\
& G_+Qe_\alpha=e_\alpha, \notag \\
& G_+QG_-e_\alpha =G_-e_\alpha. \notag
\end{align}
We see that $[G_-,G_+]=0$; $\Pi_4=[Q,G_+]$ is the projection to $H_4$ along $H_0$; 
$\Pi_0=\mathrm{Id}-\Pi_4$ is the projection to $H_0$ along $H_4$. 

Consider the integral $\int\colon H\to\C$. 
\begin{description}
\item[A6]
We require that
\begin{align}
& \int Q(a)b = (-1)^{\tilde a+1}\int aQ(b), \\
& \int G_-(a)b = (-1)^{\tilde a}\int aG_-(b), \notag \\
& \int G_+(a)b  = (-1)^{\tilde a}\int aG_+(b). \notag
\end{align}
\end{description}
These properties imply that
$\int G_-G_+(a)b=\int aG_-G_+(b)$, $\int \Pi_4(a)b=\int a\Pi_4(b)$, and
$\int \Pi_0(a)b=\int a\Pi_0(b)$.

We can define a scalar product on $H$ as $(a,b)=\int ab.$
\begin{description}
\item[A7]
This scalar product 
is non-degenerate. 
\end{description}
Using the scalar product we may turn any graded symmetric or skew-symmetric operator $A: H \to H$
into the bivector (that we denote by $[A]$), well-defined up to a sign. 

\subsection{Tensor expressions in terms of graphs}\label{sec:tensors}
Here we explain a way to encode some tensor expressions over an arbitrary vector space in terms of graphs.

Consider an arbitrary graph (we allow graphs to have leaves and we require vertices to be at least of 
degree $3$, the definition of graph that we use can be found in~\cite{Man99b}). We associate a symmetric 
$n$-form to each internal vertex of degree $n$, a symmetric bivector 
to each egde, and a vector to each leaf. Then we can substitute the tensor product of all vectors in leaves 
and bivectors in edges into the product of $n$-forms in vertices, distributing the components of tensors in 
the same way as the corresponding edges and leaves are attached to vertices in the graph. This way we get a 
number (in the case of $\Z_2$-graded vector space there is an additional sign correction, see Seection~\ref{sec:remarks} below).

Let us study an example:
\begin{equation}
\inspic{pictures.8}
\end{equation}
We assign a $5$-form $x$ to the left vertex of this graph and a $3$-form $y$ to the right vertex. 
Then the number that we get from this graph is $x(a,b,c,v,w)\cdot y(v,w,d)$.

Note that vectors, bivectors and $n$-forms used in this construction can depend on some variables. 
Then what we get is not a number, but a function.

\subsection{Usage of graphs in cyclic Hodge algebras}\label{sec:usage}

Consider a cyclic Hodge algebra $H$. There are some standard tensors over $H$, 
which we associate to vertices, edges, and leaves of graphs below. 
Here we introduce the notations for these tensors.

We always assign the form
\begin{equation}\label{eq:internal}
(a_1,\dots,a_n)\mapsto \int a_1\cdot\dots\cdot a_n
\end{equation}
to a vertex of degree $n$.

There are two bivectors that can be assigned to edges: 
$[G_-G_+]$ and $[Id]$. The last one will is used only on loops.

The vectors that we will put at leaves depend on some variables. Let $\{e_1,\dots,e_s\}$ be a
homogeneous basis of $H_0$. In particular, we assume that $e_1$ is the unit of
$H$. To each vector $e_i$ we associate formal variables $t_{d,i}$, 
$d\geq 0$, of the same parity as $e_i$. Then we will put at a leaf one of the vectors 
$E_d=\sum_{i=1}^s e_i t_{d,i}$, $d\geq 0$.

\subsubsection{Remarks}\label{sec:remarks}

There is a subtlety related to the fact that $H$ is a $\Z_2$-graded space. 
The complete definition of the tensor contraction associated to a graph includes also a sign issue.  Suppose we consider a graph of genus $g$. We can choose $g$ edges in such a way that 
the graph being cut at these edge turns into a tree. To each of these edges 
we have already assigned a bivector $[A]$ for some operator $A\colon H\to H$. 
Now we have to put the bivector $[JA]$ instead of the bivector $[A]$, where 
$J$ is an operator defined by the formula $J\colon a\mapsto (-1)^{\tilde a} a$.

Another subtlety appears when we study infinite-dimensional cyclic Hodge algebras. 
In this case, we assume that all expressions in graphs that we use do converge.

\subsection{Correlators}\label{sec:correlatorsHFT}

We are going to define the potential using correlators. Let
\begin{equation}
\langle \tau_{d_1}(v_1)\dots \tau_{d_n}(v_n) \rangle_g
\end{equation}
be the sum over graphs of genus $g$ with $n$ leaves marked by $\tau_{d_i}(v_i)$, $i=1,
\dots, n$, where $v_1,\dots,v_n$ are some vectors in $H$.
The index of each internal vertex of these graphs is $\geq 3$; we associate to it
the symmetric form~\eqref{eq:internal}.
There are two possible types of edges: arbitrary edges marked by $[G_-G_+]$ (`heavy edges') 
and loops marked by $[Id]$ (`empty loops').

Consider a vertex of such graph. Let us describe all possible half-edges adjusted to
this vertex. There are $2g$, $g\geq 0$, half-edges coming from $g$ empty loops; $m$ 
half-edges coming from heavy edges of graph, and $l$ leaves 
marked $\tau_{d_{a_1}}(v_{a_1}), \dots, \tau_{d_{a_l}}(v_{a_l})$. Then we say that 
the type of this vertex is $(g,m;d_{a_1},\dots,d_{a_l})$. We denote the type of a 
vertex $\mathfrak{v}$ by $(g(\mathfrak{v}),m(\mathfrak{v});d_{a_1(\mathfrak{v})},\dots,d_{a_{l(\mathfrak{v})}(\mathfrak{v})})$.

Consider a graph $\Gamma$ in the sum determining the correlator
\begin{equation}
\langle \tau_{d_1}(v_1)\dots \tau_{d_n}(v_n) \rangle_g
\end{equation}
We associate to $\Gamma$ a number: we contract according to the graph structure all
tensors corresponding to its vertices, edges, and leaves (for leaves, we take
vectors $v_1,\dots,v_n$). Let us denote this number by $T(\Gamma)$.

Also we weight each graph by a coefficient which is the product of two
combinatorial constants. The first factor is equal to
\begin{equation}
A(\Gamma)=\frac{\prod_{\mathfrak{v}\in \mathfrak{V}(\Gamma)}2^{g(\mathfrak{v})}g(\mathfrak{v})!}{|\mathrm{aut}(\Gamma)|}.
\end{equation}
Here $|\mathrm{aut}(\Gamma)|$ is the order of the automorphism group of the labeled
graph $\Gamma$, $\mathfrak{V}(\Gamma)$ is the set of internal vertices of $\Gamma$. In other
words,
we can label each vertex $\mathfrak{v}$ by $g(\mathfrak{v})$, delete all empty loops, and then we get a
graph with the order of the automorphism group equal to $1/A(\Gamma)$.

The second factor is equal to
\begin{equation}
P(\Gamma)=\prod_{\mathfrak{v}\in \mathfrak{V}(\Gamma)}\int_{\oM_{g(\mathfrak{v}),m(\mathfrak{v})+l(\mathfrak{v})}}
\psi_1^{d_{a_1(\mathfrak{v})}}\dots
\psi_{l(\mathfrak{v})}^{d_{a_{l(\mathfrak{v})}(\mathfrak{v})}}.
\end{equation}

So, the whole contribution of the graph $\Gamma$ to the correlator is equal to
$P(\Gamma)A(\Gamma)T(\Gamma)$. One can check that the non-trivial contribution to 
the correlator $\langle\tau_{d_1}(v_1)\dots \tau_{d_n}(v_n) \rangle_g$ is given only 
by graphs that have exactly $3g-3+n-\sum_{i=1}^n d_i$ heavy edges.

\subsection{Potential}\label{sec:potential}

We fix a cyclic Hodge algebra and consider the formal power series $Z=Z(t_{d,i})$ defined as
\begin{multline}\label{potential}
Z=
\exp\left(\sum_{g=0}^\infty \hbar ^{g-1} F_g\right)\\
=\exp\left(\sum_{g=0}^\infty \hbar^{g-1}\sum_{n} \frac{1}{n!} \sum_{d_1,\dots,d_n \in \Z_{\geq 0}}
\langle \tau_{d_1}(E_{d_1})\dots \tau_{d_n}(E_{d_n})\rangle_g\right).
\end{multline}

Losev, me, and Shneiberg proved in~\cite{LosShaShn07} that this formal power series is indeed the
potential of a CohFT, that is, its components of fixed genus satisfy string,
dilaton, and all tautological relations (particular tautological relations were
checked before in~\cite{LosSha07, Sha07, ShaShn07, Shn08}). 

Our goal in this paper is to revise this theorem completely (both the statement
and the prove) using the Givental group action on CohFT.

\subsubsection{Trivial example}
Consider the trivial cyclic Hodge algebra: $H=H_0=\langle e_1\rangle$,
$Q=G_-=0$, $\int e_1=1$.
Then $E_a=e_1\cdot t_a$, and the correlator 
$
\langle\tau_{a_1}(E_{a_1})\dots\tau_{a_n}(E_{a_n})\rangle_g
$
consists just of one graph with one vertex, $g$ empty loops, and $n$ leaves marked by $a_1,\dots,a_n$.
The explicit value of the coefficient of this graph is, by definition,
\begin{equation}\label{GWpoint}
\langle\tau_{d_1}\dots\tau_{d_n}\rangle_g:=\int_{\oM_{g,n}}
\psi_1^{d_1}\dots\psi_{n}^{d_n}.
\end{equation}
So, in the case of the trivial cyclic Hodge algebra we obtain exactly the Gromov-Witten
potential of the point.


\section{Hodge field theory revisited}\label{sec:newhodge}

In this section, we reformulate my theorem with Losev and Shnei\-berg
mentioned in Section~\ref{sec:potential}, and give a different prove of it. 

\subsection{Input}\label{sec:input}

We start with a $\Z_2$-graded commutative Frobenius algebra $H$ with a unit $e_1\in H$.
We still require the odd operators $Q$, $G_-$, and $G_+$, where the first two
operators determine on $H$ the structure of Hodge bicomplex and the third one is
responsible for the particular choice of splitting $H=H_0\oplus H_4$, where
$e_1\in H_0$. We also require that the operator $Q$ is (graded) skewsymmetric,
the operators $G_-$ and $G_+$ are (graded) symmetric.

The structure of Frobenius algebra on $H$ is equivalent to the structure of
TFT on it, see Section~\ref{sec:TFT}. We can choose some basis in the spaces $H_0$ and $H_4$,
and use it in order to write down the potential associated to this TFT. Denote
this potential by $Z^{\circ}$.

\begin{proposition}\label{prop:input}
The Leibniz rule for $Q$~\eqref{eq:leibniz} is equivalent to the equation
$QZ^{\circ}=0$. The $7$-term relation~\eqref{eq:7-term} and
$1/12$-axiom~\eqref{eq:1/12-axiom} for $G_-$ are equivalent
to $(G_-z)\hat{\ }Z^{\circ}=0$.
\end{proposition}

\subsection{Potential}

Consider the formal power series 
\begin{equation}
Z:=\exp(-(G_-G_+z)\hat{\ })Z^\circ \vert_{H_0}.
\end{equation}
 Since it is restricted to $H_0$, we can consider it
as a formal
power series in $t_{d,i}$, $d=0,1,\dots$, $i=1,\dots,s=\dim H_0$ (where the
variables a linked to the same basis as was chosen in Section~\ref{sec:usage}).

\begin{proposition}\label{prop:potential}
The formal power series $Z$ defined here coinsides with the potential defined in
Section~\ref{sec:potential}.
\end{proposition}

We also want to give an alternative proof of the following statement.

\begin{proposition}\label{prop:CohFT}
The formal power series $Z$ is the potential of a CohFT on $H_0$.
\end{proposition}

\subsection{Proofs}

In this section, we collect the proofs of
Propositions~\ref{prop:input}--\ref{prop:CohFT}

\subsubsection{Proof of the first statement of Proposition~\ref{prop:input}}
\label{sec:proof-1}
Observe, that in the
language of graphs components $F^\circ_g$, $g\geq 0$, of the potential
$Z^{\circ}=\exp\left(\sum_{g\geq 0} \hbar^{g-1} F^\circ_g \right)$ are
represented as sums of one-vertex graphs of genus $g$. There are $g$ empty loops
attached to the vertices of these graphs. The equation $QZ^\circ=0$ means that if we apply $Q$ to each of the leaves of
any particular such graph, we get $0$.

Applying $Q$ to each of the leaves of the simplest graph, we get exactly the
Leibniz rule for $Q$:
\begin{equation}
\inspic{pictures.9}+\inspic{pictures.10}+\inspic{pictures.11}=0
\end{equation}
It immediately implies that if we apply $Q$ to a genus $0$ one-vertex graphs, we
also get $0$. 

Since $Q$ is skewsymmetric, we conclude that if we put $Q$ at two
different ends of an empty edge (loop), we get $0$. This implies that we also
get $0$ when we apply $Q$ to the graphs with an arbitrary number of empty loops.

Thus we see that the condition $QZ^\circ=0$ contains the Leibniz rule for $Q$ as a special
case, on the one hand, and the Leibniz rule for $Q$ together with the fact that
$Q$ is skewsymmetric imply that $QZ^\circ=0$.

\subsubsection{Proof of the secont statement of Proposition~\ref{prop:input}}
\label{sec:proof-2}

Now we consider the condition $(G_-z)\hat{\ }Z^\circ=0$. The operator
$(G_-z)\hat{\ }$ consists of two parts. The linear part, $\psi G_-$, increases
the power of $\psi$-class on a leaf (or, more precisely, it increases the power of
$\psi$-class in the contribution of a leaf in the coefficient $A(\Gamma)$), 
and applies $G_-$ to the input at this
leaf (in particular, since $G_-e_1=0$, the first term of
Formula~\ref{eq:derivative}
disappears). The quadratic part connects two leaves without $\psi$-classes
of the same or of two different
graphs with an edge; the bivector on this edge is $-[G_-]$.

So, $(G_-z)\hat{\ }Z^\circ=0$ is equivalent to the following. Fix genus $g$ and
the number of leaves $n$. Mark leaves with the powers of $\psi$-classes
$d_1,\dots,d_n$, such that $\sum_{i=1}^n d_i=3g-4+n$. Then consider the sum of
graphs of three possible types:
\begin{enumerate}
\item
One-vertex graphs with $g$ empty loops; the degrees of $\psi$-classes on leaves
are $d_1,\dots,d_{i-1},d_i+1,d_{i+1},\dots,d_n$, and $G_-$ is applied to the
input of $i$-th leaf.
\item
One-vertex graphs with $g-1$ empty loops and $1$ loop marked by $[G_-]$.
The degrees of $\psi$-classes on leaves are $d_1,\dots,d_n$.
\item
Two-vertex graphs with $g$ empty loops distributed in an arbitrary way between
two vertices; the edge connecting two vertices is marked by $[G_-]$; the
degrees of $\psi$-classes on leaves are $d_1,\dots,d_n$, and the leaves are
distributed between two vertices in an arbitrary way.
\end{enumerate}
Of course, all graphs that we consider are stable in the sense that the index of
each vertex is at least $3$. The property $(G_-z)\hat{\ }Z^\circ$ is equivalent to the
property that this sum is equal to $0$ for all $g$, $n$, $d_1,\dots,d_n$ as
above.

The simplest equation in this system appear when $3g-4+n=0$. This happens when
either $g=0$, $n=4$, or $g=n=1$. Let us describe these two equations
pictorially:
\begin{align}
&\inspic{pictures.12}+\inspic{pictures.13}+\inspic{pictures.14}+\inspic{pictures.15}
\label{eq:7term-graph} \\
& -\inspic{pictures.16}-\inspic{pictures.17}-\inspic{pictures.18} =0 \notag 
\end{align}
and
\begin{align}
& \inspic{pictures.19}-\inspic{pictures.20} =0
\label{eq:1/12axiom-graph}
\end{align}
Recall that each graph $\Gamma$ denotes a product of three constants related to
it, $P(\Gamma)$, $A(\Gamma)$, and $T(\Gamma)$. We compute the coefficients 
$P(\Gamma)A(\Gamma)$ and rewrite these two formulas as
\begin{align}
&T\left(\inspic{pictures.12}\right)+T\left(\inspic{pictures.13}\right)+T\left(\inspic{pictures.14}\right)
+T\left(\inspic{pictures.15}\right)
\label{eq:7term-graph-1} \\
& -T\left(\inspic{pictures.16}\right)-T\left(\inspic{pictures.17}\right)-T\left(\inspic{pictures.18}\right) =0 \notag 
\end{align}
and
\begin{align}
& \frac{1}{24}T\left(\inspic{pictures.19}\right)-\frac{1}{2}T\left(\inspic{pictures.20}\right) =0
\label{eq:1/12axiom-graph-1}
\end{align}
This
allows us to see immediately, that Equations~\eqref{eq:7term-graph}
and~\eqref{eq:1/12axiom-graph} are exactly the $7$-term
relation~\eqref{eq:7-term} and $1/12$-axiom~\eqref{eq:1/12-axiom}, repectively.

Thus we proved that $(G_-z)\hat{\ }Z^\circ=0$ implies the $7$-term relation and
the $1/12$-axiom. The converse statement is proved by essentially the same
argument as the Main Lemma (Lemma~4) in~\cite{LosShaShn07}.

\subsubsection{Proof of Proposition~\ref{prop:potential}}

Consider expression $\exp(-(G_-G_+z)\hat{\ })Z^\circ \vert_{H_0}$. We
have already mentioned in Section~\ref{sec:proof-2} the meaning of both the quadratic and linear part of 
a differential operator applied to an exponential formal power series. 
So, applying $-(G_-G_+z)\hat{\ }$ we either increase the degree of $\psi$ class
on a leaf and apply $G_-G_+$ to the vector on this leaf, or we connect two leaves
with no $\psi$-classes on them with the edge marked by $[G_-G_+]$. 

This way we get the following description for $\exp(-(G_-G_+z)\hat{\ })Z^\circ$.
It is a sum of graphs with the vertices described in Section~\ref{sec:proof-1}.
There are some edges marked by $(-\psi)^i[(G_-G_+)^{i+j+1}/i!j!(i+j+1)](-\psi)^j$, and at each leaf we apply $\exp(-\psi G_-G_+)$ to its input (increase the degree of $\psi$-class and simulteneously
apply the corresponding power of the operator). This, is fact, is a very general
description that suites for any operator $\exp(-(Oz)\hat{\ })$ such that $O$ vanishes the unit. 

Since $(G_-G_+)^2=0$, we get a substantially simpler picture. First, 
$(-\psi)^i[(G_-G_+)^{i+j+1}/i!j!(i+j+1)](-\psi)^j=[G_-G_+]$, so our edges are just
marked by $[G_-G_+]$ (heavy edges from Section~\ref{sec:correlatorsHFT}), and there 
are no $\psi$-classes on these edges (remember now that in the definition of $P(\Gamma)$
in Section~\ref{sec:correlatorsHFT}, the half-edges coming from heavy edges play
the role of marked points with no $\psi$-classes). There are still $\exp(-\psi
G_-G_+)=1-\psi G_-G_+$ on the leaves, but they disappear after the restriction
to $H_0$, since $G_-G_+H_0=0$.

The only thing that we have to mention in addition is that the autmorphisms of
graphs are always arranged automatically in all expression of this type. This
completes the pictorial description of $Z=\exp(-(G_-G_+z)\hat{\ })Z^\circ
\vert_{H_0}$ that, as we see, coinsides with the descripion of graphs and their
contributions from Section~\ref{sec:correlatorsHFT}.

\subsubsection{Proof of Proposition~\ref{prop:CohFT}}

Since we know that Givental's action preserves string, dilaton, and all 
tautological equations (\cite{FabShaZvo06} and Section~\ref{sec:action-corr}), 
the problem is only with the restriction of the
potential $\exp(-(G_-G_+z)\hat{\ })Z^\circ$ of a CohFt on $H$ to $H_0$. 

Observe that 
\begin{equation}
Q\exp(-(G_-G_+z)\hat{\ })Z^\circ = 
exp(-(G_-G_+z)\hat{\ })(Q+(G_-z)\hat{\ })Z^\circ = 0.
\end{equation}
 Now it is a general
statement, that the restriction of a $Q$-closed CohFT-potential to the
subspace representing the cohomology of $Q$ gives a CohFT-potential on this
spaces.

Indeed, let $\eta$ denote the scalar product on $H$ and $\eta_0$ denote the
restriction of the scalar product to $H_0$. The difference between the bivectors
$\eta^{-1}$ and $\eta_0^{-1}$ is $Q$-exact. Consider a tautological equation (as
a relation for correlators), that might use $\eta^{-1}$. If we replace
$\eta^{-1}$ with $\eta_0^{-1}$, and use the fact that all correlators are
$Q$-closed, we still get a tautological equation, but it is valid now only modulo
$Q$-exact terms. These $Q$-exact terms vanish when we restrict it to $H_0$.
This completes the proof.

\subsection{Remark}

It is obvious that we don't use exactly the axioms described in
Sections~\ref{sec:cH-algebras} and~\ref{sec:input}. In fact, much weaker system of
axioms is to be sufficient. For example, instead of Hodge property, we can just
require that there exists an odd symmetric operator $A$ (it is $G_-G_+$ in the
original settings) such that $[Q,exp(-zA)]=zG_-$ and $A^2=0$. This way we get
the condition that $H$ splits as $H_0\oplus H_2\oplus H_4$=$H_0
\oplus
\bigoplus_{\beta} 
\langle 
e_\beta, Qe_\beta
\rangle
\oplus
\bigoplus_{\alpha} 
\langle 
e_\alpha, Qe_\alpha, G_-e_\alpha, QG_-e_\alpha
\rangle$,
where $QH_0=G_-H_0=G_-H_2=0$. This generalization was considered
in~\cite{Par07}.
Second, one can observe that we don't use the whole power of Givental's theory,
but just a small part of it.  This also gives a hint on futher generalizations
of BCOV theory, its homotopy version, possible relations to the Zwiebach-type
induction procedure considered in~\cite{LosSha07,LosShaShn07}, and so on. We are
going to discuss all these questions in the next paper.



\begin{thebibliography}{00}

\bibitem{Bar01} S.~Barannikov, Quantum periods. I. Semi-infinite variations of Hodge structures,
 Internat. Math. Res. Notices  \textbf{2001},  no. 23, 1243--1264.
 
\bibitem{Bar02} S.~Barannikov,  Non-commutative periods and mirror symmetry in higher dimensions,  
Comm. Math. Phys.  \textbf{228}  (2002),  no. 2, 281--325.

\bibitem{BarKon98} S. Barannikov, M. Kontsevich, Frobenius manifolds and formality of 
Lie algebras of polyvector fields, Internat. Math. Res. Notices \textbf{1998}, no. 4, 201-215.

\bibitem{BelPan00} P.~Belorousski, R.~Pandharipande, A descendent relation in genus 2,
Ann. Scuola Norm. Sup. Pisa Cl. Sci. (4) \textbf{29} (2000), no. 1, 171--191.

\bibitem{BerCecOogVaf94} M. Bershadsky, S. Cecotti, H. Ooguri, C. Vafa, Kodaira-Spencer theory of gravity and exact 
results for quantum string amplitudes, Comm. Math. Phys. \textbf{165} (1994), no. 2, 311-427.

\bibitem{Cos07a} K.~Costello, Renormalization of the BCOV action, talk at
Augsburg, May 2007.

\bibitem{Cos07b} K.~Costello,  Renormalisation and the Batalin-Vilkovisky
formalism, arXiv: 0706.1533.

\bibitem{Dij95} R.~Dijkgraaf, Mirror symmetry and elliptic curves,
 The moduli space of curves (Texel Island, 1994),  149--163, Progr. Math., 129, Birkhäuser Boston, Boston, MA, 1995.

\bibitem{Dub96} B.~Dubrovin, Geometry of $2$D topological field theories.  
Integrable systems and quantum groups (Montecatini Terme, 1993),  
120--348, Lecture Notes in Math., 1620, Springer, Berlin, 1996.

\bibitem{DubZha01} B.~Dubrovin, Y.~Zhang, Normal forms of hierarchies of integrable PDEs, 
Frobenius manifolds and Gromov-Witten invariants, arXiv: math.DG/0108160.

\bibitem{EguXio98}
T.~Eguchi, C.-S.~Xiong, 
Quantum cohomology at higher genus: topological recursion relations and Virasoro conditions,
Adv. Theor. Math. Phys.  \textbf{2}  (1998),  no. 1, 219--229.

\bibitem{FabPan05} C.~Faber, R.~Pandharipande, 
Relative maps and tautological classes,  
J. Eur. Math. Soc. (JEMS)  \textbf{7}  (2005),  no. 1, 13--49.

\bibitem{FabShaZvo06}
C.~Faber, S.~Shadrin, D.~Zvonkine, 
Tautological relations and the r-spin Witten conjecture, arXiv:math/0612510.

\bibitem{Get94} E. Getzler, Batalin-Vilkovisky algebras and two-dimensional 
topological field theories, Comm. Math. Phys. \textbf{159} (1994), 265-285.

\bibitem{Get97} E.~Getzler, Intersection theory on $\oM_{1,4}$ and 
elliptic Gromov-Witten invariants.  J. Amer. Math. Soc.  \textbf{10}  (1997),  no. 4, 973--998.

\bibitem{Get98} E.~Getzler, Topological recursion relations in genus $2$, 
Integrable Systems and Algebraic Geometry (Kobe/Kyoto, 1997), World Scientific, 
River Edge, NJ, 1998, pp. 73--106.

\bibitem{Get04} E.~Getzler, 
The jet-space of a Frobenius manifold and higher-genus Gromov-Witten invariants, 
Frobenius manifolds, 45--89,
Aspects Math., E36, Vieweg, Wiesbaden, 2004. 

\bibitem{GetKap98} E.~Getzler, M.~Kapranov, Modular operads,
Compositio Math. \textbf{110}  (1998),  no. 1, 65--126.


\bibitem{Giv01a}
A.~Givental, Semisimple Frobenius structures at higher genus,  Internat. Math.
Res. Notices  \textbf{2001},  no. 23, 1265--1286.

\bibitem{Giv01} 
A. Givental,
Gromov-Witten invariants and quantization of quadratic 
hamiltonians,
Mosc. Math. J.  \textbf{1}  (2001),  no.~4, 551--568.


\bibitem{Giv04} 
A. Givental,
Symplectic geometry of Frobenius structures,
In ``Frobenius manifolds'',  91--112, 
Aspects Math., E36, Vieweg, Wiesbaden, 2004.

\bibitem{Ion02} E.-N.~Ionel, Topological recursive relations in
$H^{2g}(\M_{g,n})$, Invent. Math.  \textbf{148}  (2002),  no. 3, 627--658. 

\bibitem{Kaz07a} M.~Kazarian, privite communications, 2007.

\bibitem{Kaz07}
M.~Kazarian, Deformations of cohomological field theories, preprint 2007.

\bibitem{Kee92} S.~Keel, Intersection theory of moduli space of stable 
$n$-pointed curves of genus zero.  
Trans. Amer. Math. Soc. \textbf{330} (1992),  no. 2, 545--574.

\bibitem{KimLiu06} T.~Kimura, X.~Liu, A genus-3 topological recursion relation. 
Comm. Math. Phys. \textbf{262} (2006), no. 3, 645--661.

\bibitem{KimStaVor95} T.~Kimura, J.~Stasheff, A.~Voronov, 
On operad structures of moduli spaces and string theory.  Comm. Math. Phys. 
\textbf{171}  (1995),  no. 1, 1--25. 

\bibitem{Lee03} Y.-P.~Lee,
Witten's conjecture, Virasoro conjecture, and invariance of tautological
equations, arXiv:math/0311100.

\bibitem{Lee06a} Y.-P.~Lee,
Invariance of tautological equations I: conjectures and applications, 
arXiv:math/0604318.
 
\bibitem{Lee06b} Y.-P.~Lee,
Invariance of tautological equations II: Gromov--Witten theory,
arXiv:math/0605708.

\bibitem{Lee07} Y.-P.~Lee,
Notes on axiomatic Gromov--Witten theory and applications,
arXiv:0710.4349.

\bibitem{Los98}
A. Losev, Hodge strings and elements of K. Saito's theory of primitive form.  
Topological field theory, primitive forms and related topics (Kyoto, 1996),  
305--335, Progr. Math., 160, Birkhäuser Boston, Boston, MA, 1998.

\bibitem{LosSha07} A.~Losev, S.~Shadrin, From Zwiebach invariants to Getzler relation, 
Comm. Math. Phys. \textbf{271} (2007), no. 3, 649--679.

\bibitem{LosShaShn07} A.~Losev, S.~Shadrin, I.~Shneiberg, Tautological relations in Hodge field theory, 
 Nuclear Phys. B  \textbf{786}  (2007),  no. 3, 267--296.

\bibitem{Man99a} Yu. I. Manin, Three constructions of Frobenius manifolds:a comparative study, 
Asian J. Math. \textbf{3} (1999), no. 1, 179--220.

\bibitem{Man99b} Yu.~Manin, Frobenius manifolds, quantum cohomology, and moduli spaces. 
American Mathematical Society Colloquium Publications, 47. American Mathematical Society, 
Providence, RI, 1999.

\bibitem{Mer98} S.~Merkulov, 
Formality of canonical symplectic complexes and Frobenius manifolds.  Internat.
Math. Res. Notices  \textbf{1998},  no. 14, 727--733.

\bibitem{Mum83}
D.~Mumford, Towards an enumerative geometry of the moduli space of curves in Arithmetic and geometry II 
(M. Artin, J. Tate, eds.), Progress in Math., vol. 36, Birkh\"auser, Basel, 1983.

\bibitem{Par07}
J.-S.~Park, 
Semi-classical quantum field theories and Frobenius manifolds,  Lett. Math.
Phys.  \textbf{81}  (2007),  no. 1, 41--59.

\bibitem{Sha07} S.~Shadrin, A definition of descendants at one point in graph calculus, 
 Adv. Theor. Math. Phys.  \textbf{11}  (2007),  no. 3, 351--370.

\bibitem{ShaShn07} S.~Shadrin, I.~Shneiberg, Belorousski-Pandharipande relation in dGBV algebras,
J. Geom. Phys. \textbf{57} (2007), no. 2, 597--615.

\bibitem{Shn08} I.~Shneiberg, Topological recursion relation in $\oM_{2,2}$,
Funct. Anal. Appl. \textbf{42} (2008), to appear.

\bibitem{Tel07}  C.~Teleman, The structure of 2D semi-simple field theories,
arXiv:0712.0160.

\bibitem{Wit01} E. Witten. Two dimensional gravity and intersection theory on moduli space, 
Surveys in Differential Geometry, vol.~1 (1991), 243--310.

\bibitem{Zwi93} B.~Zwiebach, Closed string field theory: quantum action and the Batalin-Vilkovisky 
master equation, Nuclear Phys. \textbf{B390} (1993), no. 1, 33--152.

\end{thebibliography}
\end{document}